\documentclass[12pt]{article}
\usepackage{amsfonts}
\usepackage{}
\usepackage{mathrsfs,graphicx,amsmath,amssymb,amsfonts}

\usepackage{appendix}
\usepackage{caption}
\usepackage{mathrsfs,graphicx,amsmath,amssymb,amsfonts,amsthm}
\usepackage{cite}
\usepackage{longtable}
\usepackage{supertabular}
\evensidemargin=\oddsidemargin\topmargin=-1.5cm

\theoremstyle{definition}

\begin{document}
\title{The distance spectrum of the complements of graphs of diameter greater than three \thanks{This work is supported by NSFC (No. 11461071).}}
\author{Xu Chen, Guoping Wang\footnote{Corresponding author. Email: xj.wgp@163.com.}\\
{\small School of Mathematical Sciences, Xinjiang Normal University,}\\
{\small Urumqi, Xinjiang 830017, P.R.China}}
\date{}
\maketitle {\bf Abstract.}
Suppose that $G$ is a connected simple graph with the vertex set
$V( G ) = \{ v_1,v_2,\cdots ,v_n  \} $.
Let $d_G( v_i,v_j ) $ be the distance between $v_i$ and $v_j$ of $G$.
Then the distance matrix of $G$ is $D( G ) =( d_{ij} ) _{n\times n}$,
where $d_{ij}=d_G( v_i,v_j ) $.
Since $D( G )$ is a non-negative real symmetric matrix,
its eigenvalues can be arranged
$\lambda_1(G)\ge \lambda_2(G)\ge \cdots \ge \lambda_n(G)$,
where eigenvalues $\lambda_1(G)$ and $\lambda_n(G)$
are called the distance spectral radius and the least distance eigenvalue of $G$, respectively.
The {\it diameter} of graph $G$ is the farthest distance between all pairs of vertices.
In this paper, we determine the unique graph whose distance spectral radius attains maximum and minimum 
among all complements of graphs of diameter greater than three, respectively.
Furthermore, we also characterize the unique graph whose least distance eigenvalue attains maximum and minimum 
among all complements of graphs of diameter greater than three, respectively.

{\flushleft{\bf CLC number:}}  O 157.5\\

\section{Introduction}

~~~~The distance spectral radius of graphs have been studied extensively.
S. Bose, M. Nath and S. Paul \cite{B.S} determined the unique graph with maximal distance
spectral radius among graphs without a pendant vertex.
A. Ilic \cite{I.A} attained the unique graph whose distance spectral radius is maximum 
among $n$-vertex trees with perfect matching and fixed maximum degree.
W. Ning, L. Ouyang and M. Lu \cite{N.W} characterized the graph with minimum distance spectral radius 
among trees with given number of pendant vertices.
For more about the distance spectra of graphs see the survey \cite{A.M} as well as the references therein.

The least distance eigenvalues of connected graphs have been also studied.
H. Q. Lin \cite{L.H.Q.1} gave an upper bound on the least distance eigenvalue and
characterized all the connected graphs with the least distance eigenvalue
in $[-1-\sqrt{2}, a]$, where $a$ is the smallest root of $x^3-x^2-11x-7=0$ and $a\in (-1-\sqrt{2},-2)$.
H. Y. Lin and B. Zhou \cite{L.H.Y} attained the trees with the least distance eigenvalues
in $ [ -3-\sqrt{5},-2-\sqrt{2}  ] $ and the unicyclic and bicyclic graphs with least distance
eigenvalues in $(-2-\sqrt{2},-2.383)$.
G. L. Yu \cite{Y.L.G} introduced all the graphs with the least distance
eigenvalue in $[-2.383,0]$.

The complement of graph $G=( V( G ) ,E( G ) )$ is denoted by $G^c=( V( G^c ) ,E( G^c ) )$,
where $V( G^c ) =V( G ) $ and $E( G^c ) = \{ xy\notin E(G):x,y\in V( G ) ,\} $.
Y. X. Fan, F. F. Zhang and Y. Wang \cite{F.Y.Z} determined the connected graph with the minimal least eigenvalue 
among all complements of trees.
G. S. Jiang, G. D. Yu, W. Sun and Z. Ruan \cite{J.G.S} gave the graph with the minimal least eigenvalue 
among all graphs whose complements are connected and have only three pendent vertices.
S. C. Li and S. J. Wang \cite{L.S.C} introduced the unique connected graph whose least signless Laplacian
eigenvalue attains the minimum in the set of the complements of all trees.
G. D. Yu, Y. Z. Fan and M. L. Ye \cite{Y.G.D} achieved the unique graph which minimizes the least signless Laplacian eigenvalue 
among all connected complements of unicyclic graphs.

Currently there is very little research about the distance eigenvalues of complements of graphs.
H. Q. Lin and S. Drury \cite{L.H.Q.2}
characterized the unique graphs whose distance spectral radius have maximum and minimum among all complements of trees,
and the unique graphs whose least distance eigenvalue have maximum and minimum among all complements of trees.
R. Qin, D. Li, Y. Y. Chen and J. X. Meng \cite{Q.R}
determined the unique graph which has maximum distance spectral radius 
among all complements of unicyclic graphs and the unique graph which has maximum least distance eigenvalue 
among all complements of unicyclic graphs of diameter three.

Let $G$ be a connected simple graph with the vertex set $V( G ) = \{ v_1,v_2,\cdots ,v_n  \} $.
Then the adjacency matrix of $G$ is $A( G ) =( a_{ij} ) _{n\times n}$,
where $a_{ij}=1$ if $v_i$ is adjacent to $v_j$, and $a_{ij}=0$ otherwise.
In this paper, we observe the relations between $D(G^c)$ and $A(G)$ and use them
to determine the unique graph whose distance spectral radius attains maximum and minimum 
among all complements of graphs of diameter greater than three, respectively.
Furthermore, we also characterize the unique graph whose least distance eigenvalue attains maximum and minimum 
among all complements of graphs of diameter greater than three, respectively.

\section {The distance spectral radius of the complements of graphs of diameter greater than three}
~~~~The below Lemma 2.1 reflects the relationship of $D(G^c)$ and $A(G)$.  \vskip 3mm

{\noindent \bf Lemma 2.1.}  {\it Suppose that $G$ is a simple graph on $n$ vertices whose diameter $d(G)$ is greater than three.
Then we have
\begin{enumerate}
	\setlength{\parskip}{0ex}
	\item[\rm (I.)]
	when $d( G ) > 3$, $D( G^c ) =J_n-I_n+A( G )$.
	\item[\rm (II.)]
	when $d( G ) =3$, $D( G^c )\geq J_n-I_n+A( G )$.
\end{enumerate} }

{\noindent \bf Proof.}
When $d(G)>3$, for any two vertices $u$ and $v$ of $G$,
there must exist the vertex $w$ of $G$ which is adjacent to neither $u$ nor $v$.
Thus $d_{G^c}(u,v)=2$ if $u$ is adjacent to $v$ in $G$,
and $d_{G^c}(u,v)=1$ otherwise.
This shows that $D(G^c)=J_n-I_n+A(G)$.

When $d(G)=3$, take two arbitrary vertices $u$ and $v$ of $G$.
Suppose that all vertices in $G\setminus \{u, v\}$ are adjacent to either $u$ or $v$.
Then $d_{G^c}(u,v)=3$ if $u$ is adjacent to $v$,
and $d_{G^c}(u_i,u_j)=1$ otherwise.
So We assume that there exists the vertex $w$ in $G\setminus \{u, v\}$ which is adjacent to neither $u$ nor $v$.
Then $d_{G^c}(u,v)=2$ if $u$ is adjacent to $v$, and $d_{G^c}(u,v)=1$ otherwise.
This shows that $D(G^c)\ge J_n-I_n+A(G)$. $\Box$ \vskip 3mm

In this section using the relations between $D(G^c)$ and $A(G)$ stated in Lemma 2.1
we determine the unique graph whose distance spectral radius attains maximum and minimum
among all complements of graphs of diameter greater than three, respectively.
\vskip 3mm

Suppose $G$ is a connected simple graph with
the vertex set $V(G) =\{ v_1,v_2,\cdots ,v_n \}$.
Let $x=( x_1,x_2,\cdots ,x_n ) ^T$ be an eigenvector of $D( G )$ with respect to the eigenvalue $\rho$,
where $x( v_i ) =x_i$ ($i=1,2,\cdots ,n$).
Then we have
\begin{equation}
  \rho x_i=\sum_{v_j\in V( G )}{d_{ij}x_j}.
\end{equation}

Suppose that $G$ is a connected simple graph.
In this paper we always assume that its complement $G^c$ is also connected.

Let the vertex $u$ connect the $s$ vertices of the complete graph $K_{n-2}$ and $v$  connect other $t$ ($=n-2-s$) vertices of $K_{n-2}$.
We denote by $H(s,t)$ the resulting graph.
\vskip 3mm

{\noindent \bf Lemma 2.2.} {\it Suppose that $G$ is a simple graph of diameter greater than three on $n$ vertices,
and let $H(s,t)$ be the graph defined above.
Then $\lambda_1(G^c) < \lambda _1(H^c(s,t))$.}
\vskip 3mm

{\noindent \bf Proof.}  Since $d( G ) >3$,
there must be two vertices $u$ and $v$ of $G$ such that $d_G( u,v ) >3$.
Clearly, the neighbours $N_G(u)$ and $N_G(v)$ of vertices $u$ and $v$ in the graph $G$ satisfy $N_G(u)\cap N_G(v)=\emptyset$. 
Set $W=V(G)\backslash (N_G(u)\cup N_G(v)\cup \{u,v\})$.
Suppose that $s$ and $t$ are two positive integers such that $s\geq |N_G(u)|$, $t\geq |N_G(v)|$ and $s+t=n-2$.
Connecting all pairs of vertices of $G$ but $u$ and $v$,
connecting $u$ with $s-|N_G(u)|$ vertices of $W$ and connecting $v$ with other $t-|N_G(v)|$ vertices of $W$.
Obviously, the resulting graph is isomorphic to $H(s,t)$.

Let $x$ be the unit Perron vector of $D( G^c )$ with respect to $\lambda_1(G^c)$.
That is, each entry of $x$ is positive and $\|x\|=1$.
Note that $d( G ) >3$ and $d( H(s,t) ) =3$. From Lemma 2.1 we have
\begin{equation}
\begin{split}
\lambda_1(G^c) &=x^TD(G^c)x\\
&=x^T( J_n-I_n ) x+x^TA( G ) x\\
&<x^T( J_n-I_n ) x+x^TA( H(s,t) ) x\\
&\le x^TD( H^c(s,t) ) x.
\nonumber
\end{split}
\end{equation}
By Rayleigh's theorem we know $\lambda _1(H^c(s,t)) \ge x^TD( H^c(s,t) ) x$.
Then
$\lambda_1(G^c) < \lambda _1(H^c(s,t))$.
$\Box$ \vskip 3mm

Suppose that two graphs $G$ and $H$ are isomorphic. Then we write $G\cong H$.\vskip 3mm

{\noindent \bf Lemma 2.3.} {\it Let $H(s,t)$ be the graph as above.
Then
$\lambda _1(H^c(s,t)) \le \lambda _1(H^c(\lfloor\frac{n}{2}-1\rfloor, \lceil\frac{n}{2}-1\rceil))$
with equality if and only if
$H(s,t)\cong H(\lfloor\frac{n}{2}-1\rfloor, \lceil\frac{n}{2}-1\rceil)$.} \vskip 3mm

{\noindent \bf Proof.}
Let $k=\lambda _1(H^c(s,t))$.
Set $x$ to be the Perron vector of $D(H^c(s,t))$ with respect to $k$.
By the symmetry of $H^c(s,t)$ all the vertices in $N_{H(s,t)}(u)$ correspond to the same value $x_1$
and all the vertices in $N_{H(s,t)}(v)$ correspond to the same value $x_2$.
Let $x(u)=x_u$ and $x(v)=x_v$.
Then from the eigen-equation (1) we have
$$ \left\{ \begin{array}{l}
	kx_u=2sx_1+tx_2+x_v,\\
    kx_1=2x_u+2( s-1 ) x_1+3tx_2+x_v,\\
    kx_2=x_u+3sx_1+2( t-1 ) x_2+x_v,\\
    kx_v=x_u+sx_1+2tx_2.
\end{array}\right. $$
We can transform the above equations into a matrix equation $( kI_4-D ) x'=0$,
where $x'=( x_u,x_1,x_2,x_v ) ^T$ and
$$D=\left( \begin{matrix}
	0&		2s&		t&		1\\
	2&		2( s-1 )&		3t&		1\\
	1&		3s&		2( t-1 )&		2\\
	1&		s&		2t&		0\\
\end{matrix} \right).$$
Let $\phi _{s,t}( \lambda ) =\det ( I_4\lambda -D ) $. Then
\begin{equation}
\begin{split}
\phi _{s,t}( \lambda ) = & \lambda ^4+( -2s-2t+4 ) \lambda ^3+( -9s-9t-5st+3 ) \lambda ^2\\
&+ ( -12s-12t-4st-4 ) \lambda +( -4s-4t-4 ).
\nonumber
\end{split}
\end{equation}

Therefore, we obtain
$\phi _{s,t}( \lambda ) -\phi _{s-1,t+1}( \lambda ) =\lambda ( s-t-1 ) ( 5\lambda +4 )$.
Since the path $P_2$ of order 2 is an induced subgraph of $H^c(s,t)$,
$D(H^c(s,t))$ contains $D(P_2)$ as a principal submatrix.
Whereas $\lambda _1(P_2)=1$,
by Interlacing theorem we attain $\lambda _1( H^c(s,t))>1$.
Without loss of generality we assume $s\le t$.
We can compute out that $\phi _{s,t}( \lambda ) -\phi _{s-1,t+1}( \lambda ) < 0$ if $\lambda>1$.
Then $\lambda _1(H^c(s,t)) > \lambda _1( H^c(s-1,t+1))$.
Note that $s+t=n-2$. We finally obtain that
$\lambda _1(H^c(s,t)) \le \lambda _1(H^c(\lfloor\frac{n}{2}-1\rfloor, \lceil\frac{n}{2}-1\rceil))$.  $\Box$ \vskip 3mm

Combining Lemmas 2.2 and 2.3 we obtain the following main result. \vskip 3mm

{\noindent \bf Theorem 2.4.} {\it
Let $G$ be a simple graph of diameter greater than three on $n$ vertices.
Then $\lambda_1(G^c)< \lambda_1(H^c(\lfloor\frac{n}{2}-1\rfloor, \lceil\frac{n}{2}-1\rceil))$.} \vskip 3mm

We denote by $G'$ the connected graph obtained from $G$ by deleting an edge of $G$ which are not adjacent. \vskip 3mm

{\noindent \bf Lemma 2.5.} {\it
Suppose that $G$ is a simple graph of diameter greater than three on $n$ vertices, and $G'$ is above.
Then $\lambda_1(G^c)\le \lambda_1(G'^c)$.} \vskip 3mm

{\noindent \bf Proof.}
Let $y$ be the unit Perron vector of $D(G'^c)$ with respect to $\lambda_1(G'^c)$.
Note that $d(G')>d(G)>3$.
From Lemma 2.1 we have
\begin{equation}
\begin{split}
\lambda_1(G'^c) &=y^TD(G'^c)y\\
&= y^T( J_n-I_n ) y+y^TA(G') y\\
&\le y^T( J_n-I_n ) y+y^TA(G) y\\
&= y^TD(G^c) y.
\nonumber
\end{split}
\end{equation}
By Rayleigh's theorem, $\lambda _1(G^c) \ge y^TD(G^c ) y$,
and so $\lambda_1(G^c)\ge \lambda _1(G'^c)$. $\Box$ \vskip 3mm

By repeatedly applying Lemma 2.5 we can prove that the result is true. \vskip 3mm

{\noindent \bf Lemma 2.6.} {\it
Suppose that $G$ is a simple graph of diameter greater than three on $n$ vertices,
and $T$ is a spanning tree of $G$.
Then $\lambda_1(G^c)\ge \lambda_1(T^c)$
with equality if and only if $G\cong T$.}
\vskip 3mm

{\noindent \bf Lemma 2.7.} (\cite{L.H.Q.2} )  {\it
Let $P_n$ be the path of order $n\ge 4$.
If $T$ is not isomorphic to the star graph $K_{1,n-1}$,
then $\lambda_1(T^c)\ge \lambda_1(P_n^c)$ with equality if and only if $T\cong P_n$.}
\vskip 3mm

Combining Lemmas 2.5, 2.6 and 2.7 we attain the following main result.
\vskip 3mm

{\noindent \bf Theorem 2.8.} {\it
Let $G$ be a simple graph of diameter greater than three on $n$ vertices.
Then $\lambda_1(G^c)\ge \lambda_1(P_n^c)$
with equality if and only if $G\cong P_n$.}

\section {The least distance eigenvalue of the complements of graphs of diameter greater than three}
~~~~In this section using the relations between $D(G^c)$ and $A(G)$ declared in Lemma 2.1
we determine the unique graph whose least distance eigenvalue attains maximum and minimum
among all complements of graphs of diameter greater than three, respectively. \vskip 3mm

Let $T(a,b)$ denote the tree obtained from the path $P_2$ of order $2$ by appending $a$ vertices to one vertex of $P_2$ and $b$ vertices to the other.
We denote by $T_1(a,b)$ the tree obtained from $P_3$ of order $3$ by appending $a$ vertices to one end vertex of $P_3$ and $b$ vertices to the other.
Let $T_2(a,b)$ be the tree obtained from $T(a,b)$ by appending an additional pendent edge to the group of $a$ pendent vertices of $T(a,b)$.
Clearly, $d(T(a+1,b))+1=d(T_1(a,b))=d(T_2(a,b))=4$. \vskip 3mm

{\noindent \bf Lemma 3.1.} (\cite{L.H.Q.2} ) {\it
Let $T(a+1,b)$, $T_1(a,b)$ and $T_2(a,b)$ be three trees of order $n$ ($=a+b+3$) as above.
Then we have
$$\lambda_n(T_1^c(a,b)\ge \lambda_n(T_2^c(a,b)))>\lambda_n(T^c(a+1,b)).$$
The equality holds if and only if $T_1(a,b)\cong T_2(a,b)$.} \vskip 3mm

Suppose that $G$ is a simple graph of diameter greater than three with
the vertex set $V(G) =\{ v_1,v_2,\cdots ,v_n\}$ ($n\ge 7$).
Let $x=\left( x_1,x_2,\cdots ,x_n \right) ^T$ be an eigenvector of $D\left( G^c \right) $ with respect to $\lambda _n(G^c)$,
where $x\left( v_i \right) =x_i$ ($i=1, 2,\cdots, n$).
Write $V_+=\left\{ v_i\in V\left( G^c \right): x_i>0 \right\} $,
$V_-=\left\{v_i\in V\left( G^c \right): x_i<0 \right\} $
and $V_0=\left\{ v_i\in V\left( G^c \right): x_i=0 \right\} $.
Let $|V_+\cup V_0|=p$ and $|V_-|=q$.
Without loss of generality in what follows we assume that $p\ge q$.
Note that $p+q=n\ge 7$. We have $p\ge 4$. \vskip 3mm

{\noindent \bf Lemma 3.2.} {\it Suppose that $G$ is a simple graph of diameter greater than three on $n\ge 7$ vertices.
If $q=1$ then $\lambda_n(G^c)\ge \lambda_n(T^c(n-3,1))$.} \vskip 3mm

{\noindent \bf Proof.}
Note that $q=1$. We let $V_-=\{v\}$.
Since $d(G)>3$, there must exist either the path $P_4=v\widetilde{u}_1\widetilde{u}_2\widetilde{u}_3$ or the path $P_5=u_1u_2vu_4u_5$.
Now we distinguish two cases as follows.

{\bf Case 1.} Suppose that there is the path $P_4=v\widetilde{u}_1\widetilde{u}_2\widetilde{u}_3$.

In this case deleting all edges in $G\setminus \{v\}$ except $\widetilde{u}_1\widetilde{u}_2$ 
and $\widetilde{u}_2\widetilde{u}_3$ and connecting all pairs of vertices 
which are not adjacent between the vertex $v$ and all vertices of $G\setminus \{v,\widetilde{u}_1,\widetilde{u}_2,\widetilde{u}_3\}$ in $G$.
Thus we obtain the resulting graph which is isomorphic to the graph $T_1(n-3,2)$.
From the above argument we know that $x^T A(G)x=\sum_{v_iv_j\in V(G)}x_ix_j\ge \sum_{v_iv_j\in V(T_1(n-4,1))}x_ix_j=x^T A(T_1(n-4,1))x$.

{\bf Case 2.} Suppose that there is the path $P_5=u_1u_2vu_4u_5$.

In this case deleting all edges in $G\setminus \{v\}$ except $u_1u_2$ and $u_4u_5$ and connecting all pairs of vertices 
which are not adjacent between the vertex $v$ and all vertices of $G\setminus \{u_1,u_2,v,u_4,u_5\}$ in $G$.
Thus we obtain the resulting graph which is isomorphic to the graph $T_2(n-4,1)$.
From the above argument we know that $x^T A(G)x=\sum_{v_iv_j\in V(G)}x_ix_j\ge \sum_{v_iv_j\in V(T_2(n-4,1))}x_ix_j=x^T A(T_2(n-4,1))x$.

Set $x$ to be the unit eigenvector of $D(G^c)$ with respect to $\lambda_n(G^c)$.
Note that $d(G)> 3$ and $d(T_1(n-3,2))=4$.
From Lemma 2.1 we have
\begin{equation}
	\begin{split}
		\lambda_n(G^c)&=x^TD(G^c)x\\
		&=x^T(J_n-I_n)x+x^TA(G)x\\
		&\ge x^T(J_n-I_n)x+x^TA(T_1(n-4,1))x\\
		&=x^TD(T_1^c(n-4,1))x.
		\nonumber
	\end{split}
\end{equation}
By Rayleigh's theorem we attain that
$\lambda_n(T_1^c(n-4,1))\le x^TD(T_1^c(n-4,1))x$.
Then we have $\lambda_n(G^c)\ge \lambda_n(T_1^c(n-4,1))$.

Similarly, we can determine that $\lambda_n(G^c)\ge \lambda_n(T_2^c(n-4,1))$.

By Lemma 3.1, $\lambda_n(T_1^c(n-4,1)\ge \lambda_n(T_2^c(n-4,1)))>\lambda_n(T^c(n-3,1))$.
From the above argument we have $\lambda_n(G^c)> \lambda_n(T^c(n-3,1))$. $\Box$ \vskip 3mm

Let $B_1(p,q)$ be the graph obtained from the complete bipartite graph $K_{p,q}$ by deleting the edge $uv$.
Suppose $u$ and $w$ are two vertices of the partition $U$ and $v$ belongs to the partition $V$.
Deleting all edges of $B_1(p,q)$ which are incident to $w$ except $wv$ we denote by $B_2(p,q)$ the resulting graph.
Clearly, $d(B_1(p,q))+1=d(B_2(p,q))=4$.

Suppose that $S$ is a subset of $V(G)$.
Then we denote by $G[S]$ the subgraph of $G$ induced by $S$. \vskip 3mm

{\noindent \bf Lemma 3.3.} {\it Suppose that $G$ is a simple graph of diameter greater than three on $n\ge 7$ vertices.
If $q\ge 2$ then we have $\lambda_n(G^c)\ge \lambda_n(B_2^c(p,q))$.} \vskip 3mm

{\noindent \bf Proof.} Set $x$ to be the unit eigenvector of $D(G^c)$ with respect to $\lambda_n(G^c)$.
Deleting all edges in $G[V_-]$ and $G[V_+\cup V_0]$ of $G$ we denote by $G'$ the resulting bipartite graph.
If $G'$ is connected then since $d(G')\geq d(G)>3$, there must be two vertices $u$ and $w$ such that $d_{G'}(u,w)=4$.
Let $P=uu_1u_2vw$ be the path between $u$ and $w$.
Then $u$ and $w$ are in the same partition, say $u$ and $w$ are both contained in $V_+\cup V_0$.
Without loss of generality assume that $x(u)\ge x(w)$.
Deleting all edges which are incident to $w$ except $vw$ and connecting all pairs of vertices
between $(V_+\cup V_0)\setminus \{w\}$ and $V_-$ which are not adjacent except $u$ and $v$ in $G'$.
Obviously, the resulting graph is isomorphic to the graph $B_2(p,q)$.
From the above construction we know that $x^T A(G)x=\sum_{v_iv_j\in E(G)}x_ix_j\ge \sum_{v_iv_j\in E(B_2(p,q))}x_ix_j=x^T A(B_2(p,q))x$.

So we can assume that $G'$ is not connected.
Since $G$ is connected, $G'$ must have one nontrivial component, that is, it contains at least one edge.
Now we distinguish two cases as follows.

{\bf Case 1.} $G'$ has at least two nontrivial components.

Suppose two edges $\tilde{u}\tilde{v}$ and $\tilde{u}'\tilde{v}'$ belong to two distinct nontrivial components.
Without loss of generality we assume that $x(\tilde{u})\ge x(\tilde{u}')\ge 0$.
Deleting all edges which are incident to $\tilde{u}'$ except $\tilde{u}'\tilde{v}'$ and connecting all pairs of vertices
between $(V_+\cup V_0)\setminus \{\tilde{u}'\}$ and $V_-$ which are not adjacent except $\tilde{u}$ and $\tilde{v}'$ in $G'$.
Obviously, the resulting graph is isomorphic to the graph $B_2(p,q)$.
From the above construction we know that $x^T A(G)x=\sum_{v_iv_j\in E(G)}x_ix_j\ge \sum_{v_iv_j\in E(B_2(p,q))}x_ix_j=x^T A(B_2(p,q))x$.

{\bf Case 2.} $G'$ has exactly one nontrivial component.

If $G'$ has exactly one isolated vertex $\bar{w}$,
then since $d(G)>3$, there must be two vertices $\bar{u}$ and $\bar{v}$ which are not adjacent in $G'$.
Without loss of generality we assume that $x(\bar{w})\ge 0$.
Connecting all pairs of vertices between $(V_+\cup V_0)\setminus \{\bar{w}\}$ and $V_-$ which are not adjacent
except $\bar{u}$ and $\bar{v}$ and connecting $\bar{w}$ and $\bar{v}$ in $G'$.
Obviously, the resulting graph is isomorphic to the graph $B_2(p,q)$.
From the above argument we know that $x^T A(G)x=\sum_{v_iv_j\in E(G)}x_ix_j\ge \sum_{v_iv_j\in E(B_2(p,q))}x_ix_j=x^T A(B_2(p,q))x$.

So we assume that $G'$ has at least two isolated vertices $\hat{w}$ and $\hat{w}'$.
Without loss of generality we assume that $x(\hat{w})\ge 0$.
Suppose $x(\hat{v})< 0$ in the edge $\hat{u}\hat{v}$.
We denote by $G''$ the graph obtained from $G'$ by connecting all pairs of vertices
between $(V_+\cup V_0)\setminus \{\hat{w}'\}$ and $V_-\setminus \{\hat{w}'\}$ which are not adjacent except $\hat{w}$ and $\hat{v}$.
Furthermore, connect $\hat{w}'$ and $\hat{v}$ if $x(\hat{w}')\ge 0$, and connect $\hat{w}'$ and $\hat{u}$ otherwise in $G''$.
Clearly, the resulting graph is isomorphic to the graph $B_2(p,q)$.
From the above construction we know that $x^T A(G)x=\sum_{v_iv_j\in E(G)}x_ix_j\ge \sum_{v_iv_j\in E(B_2(p,q))}x_ix_j=x^T A(B_2(p,q))x$.

Note that $d(G)> 3$ and $d(B_2(p,q))=4$.
From Lemma 2.1 and the above arguments we have
\begin{equation}
\begin{split}
\lambda_n(G^c)&=x^TD(G^c)x\\
&=x^T(J_n-I_n)x+x^TA(G)x\\
&\ge x^T(J_n-I_n)x+x^TA(B_2(p,q))x\\
&=x^TD(B_2^c(p,q))x.
\nonumber
\end{split}
\end{equation}
By Rayleigh's theorem we obtain that
$\lambda_n(B_2^c(p,q))\le x^TD(B_2^c(p,q))x$. Therefore, we have
$\lambda_n(G^c)\ge \lambda_n(B_2^c(p,q))$. $\Box$ \vskip 3mm

{\noindent \bf Lemma 3.4.} {\it
Let $B_2(p,q)$ and $B_1(p,q)$ be two graphs as above.
Then we have $\lambda_n(B_1^c(p,q))<\lambda_n(B_2^c(p,q))<-3$.
} \vskip 3mm

{\noindent \bf Proof.}
Let $\lambda_n$ be the least eigenvalue of $D(B_2^c(p,q))$.
Set $x$ to be an eigenvector of $D(B_2^c(p,q))$ with respect to $\lambda_n$.
By the symmetry of $B_2^c(p,q)$ all vertices in $\{V_+\cup V_0\}\setminus \{u,w\}$
correspond to the same value $x_1$ and
all the vertices in $V_-\setminus \{v\}$ correspond to the same value $x_2$.
Set $w$ to be the only one vertex of $N_T(v)\setminus \{u\}$.
Let $x(u)=x_u$, $x(v)=x_v$ and $x(w)=x_w$.
Then from the eigen-equation (1) we obtain
$$ \left\{ \begin{array}{l}
	\lambda_nx_u=x_v+x_w+(p-2)x_1+2(q-1)x_2,\\
    \lambda_nx_v=x_u+2x_w+2(p-2)x_1+(q-1)x_2,\\
    \lambda_nx_w=x_u+2x_v+(p-2)x_1+(q-1)x_2,\\
    \lambda_nx_1=x_u+2x_v+x_w+(p-3)x_1+2(q-1)x_2,\\
    \lambda_nx_2=2x_u+x_v+x_w+2(p-2)x_1+(q-1)x_2.
\end{array}\right.$$
We can transform the above equation into a matrix equation
$(\lambda_nI_5-D_{B_2^c})x'=0$, where $x'=(x_u, x_v, x_w, x_1, x_2)$ and
$$D_{B_2^c}=\left( \begin{matrix}
	0&1&1&p-2&2(q-1)\\
	1&0&2&2(p-2)&q-1\\
	1&2&0&p-2&q-1\\
	1&2&1&p-3&2(q-1)\\
	2&1&1&2(p-2)&q-2	
\end{matrix} \right).$$
Let $\varphi_{p,q}(\lambda)=det(I_5\lambda-D_{B_2^c})$.
Then we get
\begin{equation}
\begin{split}
\varphi_{p,q}(\lambda)=
&{\lambda}^{5}- \left( q-5+p \right) {\lambda}^{4}\\
&- \left( 3\,pq+4\,p+q-10 \right) {\lambda}^{3}\\
&- \left( 8\,pq+6\,p-4\,q-8 \right) {\lambda}^{2}\\
&- \left( pq+10\,p-8 \right) \lambda+3\,pq-6\,p-2\,q+4.
\nonumber
\end{split}
\end{equation}
Similarly we have
\begin{equation}
\begin{split}
\phi_{p,q}(\lambda)=&det(I_4\lambda-D_{B_1^c})\\
=&{\lambda}^{4}+ \left( -q+4-p \right) {\lambda}^{3}\\
&+ \left( -8\,pq+2\,p+2\,q+4 \right) {\lambda}^{2}\\
&+ \left( -14\,pq+6\,p+6\,q \right) \lambda-5\,pq+2\,p+2\,q.
\end{split}
\end{equation}
By the above two equations we get
\begin{equation}
\begin{split}
\varphi_{p,q}(\lambda)-(\lambda+1)\phi_{p,q}(\lambda)=
&\left( 5\,pq-5\,p-2\,q+2 \right) {\lambda}^{3}\\
&+ \left( 14\,pq-14\,p-4\,q+4 \right) {\lambda}^{2}\\
&+ \left( 18\,pq-18\,p-8\,q+8 \right) \lambda+8\,pq-8\,p-4\,q+4.
\nonumber
\end{split}
\end{equation}

Since the path $P_4$ of order 4 is an induced subgraph of $B_2^c(p,q)$ and $B_1^c(p,q)$,
$D(P_4)$ is a principal submatrix of $D(B_2^c(p,q))$ and $D(B_1^c(p,q))$.
Whereas $\lambda_4(P_4)<-3$,
by Interlacing theorem we attain $\lambda_n(B_2^c(p,q))<-3$ and $\lambda_n(B_1^c(p,q))<-3$.
Note that $p\ge 4$ and $q\ge 2$.
We can compute out that
$\varphi_{p,q}(\lambda)-(\lambda+1)^2\phi_{p,q}(\lambda)>0$
when $\lambda<-3$.
This implies that $\lambda_n(B_2^c(p,q))>\lambda_n(B_1^c(p,q))$.
$\Box$ \vskip 3mm

{\noindent \bf Lemma 3.5.} {\it Let $B_1(p,q)$ and $T(n-3,1)$ be two graphs of order $n$ ($=p+q$) as above.
Then we have $\lambda_n(B_1^c(p,q))<\lambda_n(T^c(n-3,1))<-3$.} \vskip 3mm

{\noindent \bf Proof.}
Let $\lambda_n$ be the least eigenvalue of $D(T^c(n-3,1))$.
Set $x$ to be an eigenvector of $D(T^c(n-3,1))$ with respect to $\lambda_n$.
By the symmetry of $T^c(n-3,1)$ all vertices in $N_T(u)\setminus \{v\}$
correspond to the same value $x_1$.
Let $x(u)=x_u$, $x(v)=x_v$ and $x(w)=x_w$.
Then from the eigen-equation (1) we obtain
$$ \left\{ \begin{array}{l}
	\lambda_nx_u=3x_v+2(n-3)x_1+x_w,\\
	\lambda_nx_v=3x_u+(n-3)x_1+2x_w,\\
	\lambda_nx_1=2x_u+x_v+(n-4)x_1+x_w,\\
	\lambda_nx_w=x_u+2x_v+(n-3)x_1.
\end{array}\right.$$
We can transform the above equation into a matrix equation
$(\lambda_nI_4-D_{T^c})x'=0$, where $x'=(x_u, x_v, x_1, x_w)$ and
$$D_{T^c}=\left( \begin{matrix}
	0&3&2(n-3)&1\\
	3&0&n-3&2\\
	2&1&n-4&1\\
	1&2&n-3&0
\end{matrix} \right).$$
Let $\psi_{p,q}(\lambda)=det(I_4\lambda-D_{T^c})$.
Then we get
\begin{equation}
\begin{split}
\psi(\lambda)=
{\lambda}^{4}+ \left( -n+4 \right) {\lambda}^{3}+ \left( 4-6\,n
 \right) {\lambda}^{2}+ \left( -6\,n-8 \right) \lambda-12.
\end{split}
\end{equation}
Note that $n=p+q$. From the equations (2) and (3) we obtain
\begin{equation}
\begin{split}
		\phi_{p,q}(\lambda)-\psi(\lambda)=
		&\left( -8\,pq+8\,p+8\,q \right) {\lambda}^{2}\\
		&+\left( -14\,pq+12\,p+12\,q+8 \right) \lambda-5\,pq+2\,p+2\,q+12.
\nonumber	
\end{split}
\end{equation}

Since the Path $P_4$ of order $4$ is an induced subgraph of $T^c(n-3,1)$,
$D(T^c(n-3,1))$ contains $D(P_4)$ as a principal submatrix.
Whereas $\lambda_4(P_4)<-3$, we have $\lambda_n(T^c(n-3,1))<-3$.
Recall that $p\ge 4$ and $q\ge 2$.
Therefore, we can compute out that $\phi_{p,q}(\lambda)-\psi(\lambda)<0$ if $\lambda<-3$.
Thus, by Lemma 3.4 we get $\lambda_n(T^c(n-3,1))>\lambda_n(B_1^c(p,q))$. $\Box$ \vskip 3mm

{\noindent \bf Lemma 3.6.} {\it
Let $B_1(p,q)$ be the graph as above.
Then we have
$$\lambda_n(B_1^c(p,q))\ge \lambda_n\left(B_1^c\left(\left\lceil \frac{n}{2}\right\rceil,\left\lfloor \frac{n}{2}\right\rfloor\right)\right).$$
The equality holds if and only if
$B_1(p,q)\cong B_1(\lceil\frac{n}{2}\rceil,\lfloor \frac{n}{2}\rfloor)$.}\vskip 3mm

{\noindent \bf Proof.}
By the above equation (2) we obtain
$$\phi_{p,q}(\lambda)-\phi_{p-1,q+1}(\lambda)=\left( 8\,p-8\,q-8 \right) {\lambda}^{2}+ \left( 14\,p-14\,q-14
\right) \lambda+5\,p-5\,q-5.$$
Without loss of generality we assume that $p> q$. 
By computation we obtain that $\phi_{p,q}(\lambda)-\phi_{p-1,q+1}(\lambda)\ge 0$ if $\lambda<-3$.
Thus, by Lemma 3.4 we have $\lambda_n(B_1^c(p,q))> \lambda_n(B_1^c(\lceil \frac{n}{2}\rceil,\lfloor \frac{n}{2}\rfloor))$.
$\Box$ \vskip 3mm

Combining Lemmas 3.2-3.6 we have the following main result. \vskip 3mm

{\noindent \bf Theorem 3.7.} {\it
Let $G$ be a simple graph of diameter greater than three on $n\ge 7$ vertices.
Then we have $$\lambda_n(G^c)>
\lambda_n\left(B_1^c\left(\left\lceil \frac{n}{2}\right\rceil,\left\lfloor \frac{n}{2}\right\rfloor\right)\right).$$} \vskip 3mm


Let $x$ be the unit eigenvector of $D(G^c)$ with respect to $\lambda_n(G^c)$.
Let $G'$ denote the connected graph obtained from $G$ by deleting an edge in $V_+\cup V_0$ or $V_-$
or connecting one pair of vertices between $V_+\cup V_0$ and $V_-$ which are not adjacent such that $d(G')>3$.
Clearly, $x^TA(G)x\geq x^TA(G')x$.\vskip 3mm

{\noindent \bf Lemma 3.8.} {\it Suppose that $G$ is a simple graph of diameter greater than three on $n\ge 7$ vertices.
Then $\lambda_n(G^c)\geq \lambda_n(G'^c)$.} \vskip 3mm

{\noindent \bf Proof.}
Let $x$ be the unit eigenvector of $D(G^c)$ with respect to $\lambda_n(G^c)$.
Note that $d(G)>3$ and $d(G')>3$.
From Lemma 2.1 we have
\begin{equation}
	\begin{split}
		\lambda_n(G^c)&=x^TD(G^c)x\\
		&= x^T(J_n-I_n)x+x^TA(G)x\\
		&\ge x^T(J_n-I_n)x+x^TA(G')x\\
		&=x^TD(G'^c)x.
		\nonumber
	\end{split}
\end{equation}
By Rayleigh's theorem,
$\lambda_n(G'^c)\le x^TD(G'^c)x$, and so $\lambda_n(G^c)\ge \lambda_n(G'^c)$. $\Box$ \vskip 3mm

From Lemma 3.8 we obtained that $\lambda_n(G^c)\le \lambda_n(\check{G}^c)$ if $\check{G}$ is obtained from $G$ by
connecting one pair of vertices in $V_+\cup V_0$ or $V_-$ which are not adjacent
or deleting an edge between $V_+\cup V_0$ and $V_-$. 
Clearly, $x^TA(G)x\leq x^TA(\check{G})x$.\vskip 3mm

Let $K_{n-2}$ be a complete graph of order $n-2$.
We denote by $L'$ the graph by deleting an edge $wu$ of $K_{n-2}$ and appending a vertex $w'$ to $w$ and a vertex $v$ to $u$.
Clearly, $d(L')=4$.
Let $L''$ donote the graph by deleting an edge $w'u'$ of $K_{n-2}$ and appending a path of order $2$ to $u'$.
\vskip 3mm

{\noindent \bf Lemma 3.9.} {\it Suppose that $G$ is a simple graph of diameter greater than three on $n\ge 7$ vertices.
If $q=1$ then $\lambda_n(G^c)\le \lambda_n(L'^c)$.} \vskip 3mm

{\noindent \bf Proof.} Since $d(G)>3$, there must be two vertices $u_1$ and $u_5$ such that $d_G(u_1,u_5)=4$. 
Let $V_-=\{v\}$ and $P=u_1u_2u_3u_4u_5$. 
Now we distinguish four cases as follows.

{\bf Case 1.} Suppose that $v=u_3$ in $G$.

Without loss of generality we assume that $x(u_1)\ge x(u_5)$.
We denote by $G'$ the graph obtained from $G$ by deleting all edges which are incident to $v$ except $vu_2$,
deleting all edges which are incident to $u_5$ except $u_5u_4$ and appending them to $u_1$.
Connecting all pairs of vertices of $V(G)\setminus \{v,u_5\}$ which are not adjacent except $u_2$ and $u_4$ in $G'$.
Thus, the resulting graph is isomorphic to the graph $L'$.

{\bf Case 2.} Suppose that $v=u_2$ in $G$.

If $x(u_1)\ge x(u_5)$, we denote by $\widetilde{G}$ the graph obtained from $G$ by deleting all edges which are incident to $v$ except $vu_1$, 
deleting all edges which are indicent to $u_5$ except $u_5u_4$ and appending them to $u_1$.
Connecting all pairs of vertices of $V(G')\setminus \{v,u_5\}$ which are not adjacent except $u_1$ and $u_4$ in $\widetilde{G}$.
Obviously, the resulting graph is isomorphic to the graph $L'$.

So we assume that $x(u_1)<x(u_5)$.
We denote by $\widetilde{G}'$ the graph obtained from $G$ by deleting all edges which are incident to $v$ except $vu_3$, 
deleting all edges which are incident to $u_1$ and appending them to $u_5$.
Connecting all pairs of vertices of $V(\widetilde{G}')\setminus \{v,u_1\}$ and connecting $u_1$ and $u_5$.
Clearly, the resulting graph is isomorphic to the graph $B_2(p,q)$.

{\bf Case 3.} Suppose that $v=u_1$ in $G$.

If $x(u_2)\ge x(u_5)$, we denote by $\overline{G}$ obtained from $G$ by deleting edges which are incident to $v$ 
except $vu_2$ and deleting all edges which are incident to $u_5$ except $u_5u_4$.
Connecting all pairs of vertices of $V(\overline{G})\setminus \{v,u_5\}$ which are not adjacent except $u_2$ and $u_4$ in $\overline{G}$.
Clearly, the resulting graph is isomorphic to the graph $B_2(p,q)$.

So we assume that $x(u_2)<x(u_5)$.
We denote by $\overline{G}'$ obtained from $G$ by deleting edges which are incident to $v$ except $vu_2$ 
and deleting all edges which are incident to $u_2$ except $u_2u_3$ and $u_2v$ and appending them to $u_5$.
Connecting all pairs of vertices of $V(\overline{G}')\setminus \{v,u_2\}$ which are not adjacent except $u_3$ and $u_5$ in $\overline{G}'$.
Then the resulting graph is isomorphic to the graph $L''$.

{\bf Case 4.} Suppose that $v$ is adjacent to $u_3$ in $G$.

Without loss of generality assume that $x(u_1)\ge x(u_5)$.
We denote by $\hat{G}$ obtained from $G$ by deleting edges which are incident to $v$ except $vu_3$,
deleting all edges which are incident to $u_5$ and appending them to $u_1$ and connecting $u_5$ and $u_1$.
Connecting all pairs of vertices of $V(\hat{G})\setminus \{v,u_5\}$ except $u_1$ and $u_3$ in $\hat{G}$.
Thus, the resulting graph is isomorphic to the graph $B_2(p,q)$.

By the above arguments we obtain that the following facts.
In other cases, by deleting some edges of $G$ which are incident to $v$ we can pick out a path is isomorphic to the above four constructions.

By repeatedly applying Lemma 3.8 we can verify that
$\lambda_n(G^c)\le \lambda_n(L'^c)$ or $\lambda_n(G^c)\le \lambda_n(L''^c)$ with equality if and only if $G\cong L'$ or $G\cong L''$.

Let $\lambda_n$ be the least eigenvalue of $D(L'^c)$.
Set $x$ to be the eigenvector of $D(L'^c)$ with respect to $\lambda_n$.
By the symmetry of $L'^c$ all the vertices in $(V_+\cup V_0)\setminus \{u,v,w,w'\}$
correspond to the same value $x_1$.
Let $x(u)=x_{u}$, $x(v)=x_{v}$, $x(w)=x_w$ and $x(w')=x_{w'}$.
By the equation (1) we have
$$ \left\{ \begin{array}{l}
	\lambda_nx_u=2x_v+x_w+x_{w'}+2(p-4)x_1,\\
    \lambda_nx_v=2x_u+x_w+x_{w'}+(p-4)x_1,\\
    \lambda_nx_w=x_u+x_v+2x_{w'}+2(p-4)x_1,\\
    \lambda_nx_{w'}=x_u+x_v+2x_w+(p-4)x_1,\\
    \lambda_nx_1=2x_u+x_v+2x_w+x_{w'}+2(p-5)x_1.
\end{array}\right.$$
We can transform the above equation into a matrix equation
$(\lambda_nI_5-D_{L'^c})x'=0$, where $x'=(x_u, x_v,x_w, x_{w'}, x_1$) and
$$D_{L'^c}=\left( \begin{matrix}
	0&2&1&1&2(p-4)\\
    2&0&1&1&p-4\\
    1&1&0&2&2(p-4)\\
    1&1&2&0&p-4\\
    2&1&2&1&2(p-5)
\end{matrix} \right).$$
Let $\Psi(\lambda)=det(I_5\lambda-D_{L'^c})$. Then we get
\begin{equation}
\begin{split}
\Psi(\lambda)=det(I_5-D_{L'^c})
=&{\lambda}^{5}- \left( 2\,n-10 \right) {\lambda}^{4}\\
&- \left( -28+10\,n\right) {\lambda}^{3}-10\,n{\lambda}^{2}- \left( -4\,n+48 \right)\lambda.
\end{split}
\end{equation}
Similarly, we have
\begin{equation}
\begin{split}
\Psi'(\lambda)=det(I_5-D_{L''^c})
=&{\lambda}^{5}- \left( 2\,n-10 \right) {\lambda}^{4}- \left( -28+10\,n
 \right) {\lambda}^{3}\\
&- \left( 10+8\,n \right) {\lambda}^{2}- \left(
103-15\,n \right) \lambda+14\,n-70.
\nonumber
\end{split}
\end{equation}
From the above equations we get
$$\Psi(\lambda)-\Psi'(\lambda)=\left( -2\,n+10 \right) {\lambda}^{2}+ \left( -11\,n+55 \right)
\lambda-14\,n+70.$$

Since the path $P_5$ of order $5$ is an induced subgraph of $L'$ and $L''$, 
$D(L'^c)$ and $D(L''^c)$ contain a principal submatrix $D(P_5)$.
Whereas $\lambda_5(P_5)<-5$, by Interlacing theorem we have $\lambda_n(L'^c)<-5$ and $\lambda_n(L''^c)<-5$.
Therefore, we can compute out that $\Psi(\lambda)-\Psi'(\lambda)<0$ if $\lambda<-5$ and $n\ge 7$.
This implies that $\lambda_n(L''^c)\le \lambda_n(L'^c)$.

Thus, by the above arguments we know $\lambda_n(G^c)\le \lambda_n(L'^c)$.
$\Box$ \vskip 3mm

Let $K_p$ and $K_q$ be two complete graphs of order $p$ and $q$, respectively.
We denote by $L(p,q)$ the graph by deleting an edge $wu$ of $K_p$ and connecting $u$ and $v$ of $K_q$.
Clearly, $d(L(p,q))=4$. \vskip 3mm

{\noindent \bf Lemma 3.10.} {\it Suppose that $G$ is a simple graph of diameter greater than three on $n\ge 7$ vertices.
If $q\ge 2$ then we have $\lambda_n(G^c)\le \lambda_n(L^c(p,q))$.} \vskip 3mm

{\noindent \bf Proof.}
Note that $G$ is a connected graph.
There must be two vertices $u$ of $V_+\cup V_0$ and $v$ of $V_-$ which are adjacent in $G$.
Since $d(G)>3$, we without loss of generality assume that there exists a vertex $w$ of $V_+\cup V_0$ which are not adjacent to $u$ in $G$.
Connecting all pairs of vertices in $V_+\cup V_0$ and $V_-$ which are not adjacent 
except $w$ and $u$ and deleting all edges between $V_+\cup V_0$ and $V_-$ except $uv$ in $G$.
Obviously, the resulting graph is isomorphic to the graph $L(p,q)$.

By repeatedly applying Lemma 3.8 we can verify that the result is true. $\Box$ \vskip 3mm

{\noindent \bf Lemma 3.11.} {\it Let $L'$ and $L(p,q)$ be two graphs as above.
Then we have $\lambda_n(L'^c)<\lambda_n(L^c(p,q))$.} \vskip 3mm
	
{\noindent \bf Proof.}
Let $\lambda_n$ be the least eigenvalue of $D(L^c(p,q))$.
Set $x$ to be the eigenvector of $D(L^c(p,q))$ with respect to $\lambda_n$.
By the symmetry of $L^c(p,q)$ all the vertices in $(V_+\cup V_0)\setminus \{u,w\}$
correspond to the same value $x_1$ and
all the vertices in $V_-\setminus \{v\}$ correspond to the same value $x_2$.
Let $x(u)=x_{u}$, $x(v)=x_{v}$ and $x(w)=x_w$.
By the equation (1) we have
$$ \left\{ \begin{array}{l}
	\lambda_nx_u=2x_v+x_w+2(p-2)x_1+(q-1)x_2,\\
	\lambda_nx_v=2x_u+x_w+(p-2)x_1+2(q-1)x_2,\\
	\lambda_nx_w=x_u+x_v+2(p-2)x_1+(q-1)x_2,\\
	\lambda_nx_1=2x_u+x_v+2x_w+2(p-3)x_1+(q-1)x_2,\\
	\lambda_nx_2=x_u+2x_v+x_w+(p-2)x_1+2(q-2)x_2.
\end{array}\right.$$
We can transform the above equation into a matrix equation
$(\lambda_nI_5-D_{L^c})x'=0$, where $x'=(x_u, x_v,x_w, x_1, x_2)$ and
$$D_{L^c}=\left( \begin{matrix}
	0&2&1&2(p-2)&q-1\\
	2&0&1&p-2&2(q-1)\\
	1&1&0&2(p-2)&q-1\\
	2&1&2&2(p-3)&q-1\\
	1&2&1&p-2&2(q-2)
\end{matrix} \right).$$
Let $\Phi_{p,q}(\lambda)=det(I_5\lambda-D_{L^c})$. Then we get
\begin{equation}
	\begin{split}
		\Phi_{p,q}(\lambda)=
		&{\lambda}^{5}- \left( 2\,q-10+2\,p \right) {\lambda}^{4}\\
		&- \left( -3\,pq+16\,p+16\,q-40 \right) {\lambda}^{3}\\
		&- \left( -18\,pq+44\,p+50\,q-74\right) {\lambda}^{2}\\
		&- \left( -30\,pq+45\,p+63\,q-53 \right) \lambda+12\,pq-10\,p-22\,q+2.
	\end{split}
\end{equation}
From the equations (4) and (5) we have $$\Phi_{p,q}(\lambda)-\Psi(\lambda)= \left( 2\,p-6 \right) {\lambda}^{2}+ \left( 11\,p-33 \right)
\lambda+14\,p-42.$$

Since the path $P_5$ of order $5$ is an induced subgraph of $L^c(p,q)$,
$D(L^(p,q))$ contains $D(P_5)$ as a principal submatrix.
Whereas $\lambda_5(P_5)<-5$, by Interlacing theorem we get $\lambda_n(L^c(p,q))<-5$.
Recall that $p\ge 4$ and $q\ge 2$.
Therefore, we can compute out that $\Phi_{p,q}(\lambda)-\Psi(\lambda)> 0$ if $\lambda<-5$.
Thus, by Lemma 3.9 we obtain that $\lambda_n(L'^c)<\lambda_n(L^c(p,q))$. $\Box$ \vskip 3mm

{\noindent \bf Lemma 3.12.} {\it Let $L(p,q)$ be the graph as above.
Then we have $\lambda_n(L^c(p,q))\le \lambda_n(L^c(\lceil \frac{n}{2}\rceil,\lfloor\frac{n}{2}\rfloor))$.} \vskip 3mm

{\noindent \bf Proof.}
Note that $n=p+q$. By the equation (5) we have
\begin{equation}
	\begin{split}
\Phi_{p,q}(\lambda)-\Phi_{p-1,q+1}(\lambda)=
&\left( 3\,pq-6\,p-6\,q+12 \right) {\lambda}^{3}+ \left( 18\,pq-34\,p-
40\,q+74 \right) {\lambda}^{2}\\
&+ \left( 30\,pq-49\,p-67\,q+101 \right)
\lambda+12\,pq-10\,p-22\,q+2.
\nonumber
\end{split}
\end{equation}
Recall that $p>q\ge 2$. By computation we obtain that $\Phi_{p,q}(\lambda)-\Phi_{p-1,q+1}(\lambda)>0$ if $\lambda< -5$.

Thus, by Lemma 3.11 we have $\lambda_n(L^c(p,q))< \lambda_n(L^c(\lceil \frac{n}{2}\rceil,\lfloor\frac{n}{2}\rfloor))$. $\Box$ \vskip 3mm

Combining Lemmas 3.8-3.12 we obtain the following main result.
\vskip 3mm

{\noindent \bf Theorem 3.13.} {\it
Suppose that $G$ is a simple graph of diameter greater than three on $n\ge 7$ vertices.
Then we have
$$\lambda_n(G^c)\le \lambda_n\left(L^c\left(\left\lceil \frac{n}{2}\right\rceil,\left\lfloor\frac{n}{2}\right\rfloor\right)\right).$$}

\end{document}